\documentclass[10pt]{article}

\usepackage{latexsym}
\usepackage{amssymb}
\usepackage[frame,cmtip,arrow,matrix,line,graph]{xy}
\usepackage{epsfig}

\usepackage{amsmath}
\usepackage{amsthm}
\usepackage{amsfonts}
\usepackage{amscd}
\usepackage{color}

\newtheorem{thm}{Theorem}[section]
\newtheorem{lem}[thm]{Lemma}
\newtheorem{prop}[thm]{Proposition}

\newtheorem{claim}{Claim}

\newcommand{\vs}{\vspace{3mm}}

\newcommand{\A}{\mathcal{A}}
\newcommand{\D}{\mathcal{D}}
\newcommand{\tD}{\tilde{\mathcal{D}}}
\newcommand{\e}{\mathcal{E}}
\newcommand{\F}{\mathcal{F}}
\newcommand{\G}{\mathcal{G}}
\newcommand{\M}{\mathcal{M}}
\newcommand{\N}{\mathbb{N}}
\newcommand{\p}{\mathcal{P}}

\newcommand{\RR}{\mathbb{R}}
\newcommand{\ES}{\mathcal{S}}
\newcommand{\Z}{\mathbb{Z}}
\newcommand{\uF}{\underline{F}}
\newcommand{\uG}{\underline{G}}

\newcommand{\Top}{\textsf{Top}}

\newcommand{\Dif}{\textrm{Diff}^+}

\newcommand{\al}{\alpha}
\newcommand{\be}{\beta}
\newcommand{\ga}{\gamma}
\newcommand{\Ga}{\Gamma}

\newcommand{\de}{\delta}
\newcommand{\dop}{\Delta^{op}}
\newcommand{\eps}{\epsilon}
\newcommand{\epb}{\overline{\epsilon}}
\newcommand{\te}{\theta}

\newcommand{\la}{\lambda}

\newcommand{\Om}{\Omega}
\newcommand{\s}{\sigma}
\newcommand{\Si}{\Sigma}

\newcommand{\lar}{\longleftarrow}
\newcommand{\rar}{\longrightarrow}

\newcommand{\sta}{\stackrel}

\newcommand{\x}{\times}

\newcommand{\sq}{\square}
\newcommand{\cqfd}{\hfill $\sq$}
\newcommand{\w}{\wedge}
 
\newcommand{\bu}{\bullet} 

\newcommand{\del}{\partial}

\newcounter{samcounter}

\begin{document}

\title{Infinite loop space structure(s) on the stable mapping class group}

\author{Nathalie Wahl\\
{\small Northwestern University, 2033 Sheridan Road, Evanston IL 60208, USA}\\
{\small wahl@math.northwestern.edu}}

\maketitle

\begin{abstract}

Tillmann introduced two infinite loop space structures on the plus
 construction 
of the classifying space of the stable mapping class group,  
each with different computational advantages \cite{T1,T}. 
The first one uses disjoint union on a  
suitable cobordism category, whereas the second uses an operad  
which extends the pair of pants multiplication (i.e. the double loop  
space structure introduced by E. Y. Miller \cite{M}).  
She conjectured that these two infinite loop space structures were  
equivalent, and managed to prove that the first delooping are the  
same. In this paper, we resolve the conjecture by proving  
that the two structures are indeed equivalent, exhibiting an  
explicit geometric map.  

\end{abstract}

\section{Introduction}

Let $\Ga_{g,n}$ be the mapping class group of a  
surface $F$ of genus $g$ with $n$ boundary components.  
The classifying space $B\Ga_{g,n}$ has the homotopy type of the 
moduli space of Riemann surface of type $F$ when $n>0$. 
Attaching a torus with two boundary components to the surface  
induces  a homomorphism $\Ga_{g,1}\to \Ga_{g+1,1}$.  
Let $\Ga_\infty=\textrm{lim}_{g\to\infty}\Ga_{g,1}$ denote the  
stable mapping class group. 

 The space $\Z\x B\Ga^+_\infty$, the group completion of   
 $\coprod_{g\ge 0} B\Ga_{g,1}$, has a natural double loop space 
structure induced by the pair of pants multiplication on   
$\coprod_{g\ge 0} B\Ga_{g,1}$ (\cite{M}, see also \cite{B}).  
In \cite{T} Tillmann constructed an infinite loop space operad  
extending the pair of pants multiplication and acting on  
 $\coprod_{g\ge 0} B\Ga_{g,1}$, thus showing  
that the pair of pants multiplication actually induces an infinite loop  
space structure on $\Z\x B\Ga_\infty^+$. This multiplication 
plays a role  in conformal field theory \cite{CFT}, 
 and has also proven useful for e.g.\  
constructing homology operations \cite{B,CT}. 
Previously \cite{T1} she had 
exhibited $\Z \x B\Ga_\infty^+$ as an infinite loop space 
in quite a different way, by constructing a cobordism category $\ES$, 
symmetric  
monoidal under disjoint union of surfaces, 
such that $\Om B\ES \simeq \Z \x B\Ga_\infty^+$. 
Note that the multiplication inducing the  
infinite loop space structure in this case  
is defined on $B\ES$, and hence on a first deloop of $\Z\x B\Ga^+_\infty$. 
This infinite loop space structure has also proven useful. 
Madsen-Tillmann \cite{MT} have constructed an infinite loop map  
from $\Z\x B\Ga_\infty^+$, with the disjoint union infinite loop space 
structure, 
to $\Om^\infty \mathbb{C}P^\infty_{-1}$.   
This map has lead recently to a proof of the Mumford conjectured 
(announced by I. Madsen and M. Weiss).   
 
Tillmann conjectured \cite{T1,CT,T} that two infinite loop space  
structures were equivalent  
and managed to prove in \cite{T1} that their first 
deloopings are homotopy equivalent spaces.  
 In this paper, we resolve the conjecture of Tillmann by proving  
that the two structures are indeed equivalent, exhibiting an  
explicit geometric map. Our map sends, up to homotopy, the pair of 
pants multiplication to the loop on disjoint union multiplication 
We show that this map preserves all higher homotopies,  
using the machinery of Dwyer and Kan \cite{DK}, and hence produce an  
infinite loop map, which gives the equivalence.

\vs 

To describe our result in more details, we have to introduce some 
notation. 
Let $\M$ denote Tillmann's operad 
 and let $M$ be the associated monad.  
We will describe this operad in detail in Section \ref{opmcg}.  
The operad has $n$th space 
 $\M(n)\simeq\coprod_{g\ge 0}B\Ga_{g,n+1}$.   
So $M(*)=\M(0)\simeq\coprod_{g\ge 0}B\Ga_{g,1}$ and it is an $\M$-algebra. 
The cobordism category  $\ES$, described in Section \ref{disjoint}, has 
 objects the natural numbers.   
The morphism space 
$\ES(n,m)$ is the classifying space of a category with 
objects disjoint union of surfaces with a total of $n$  incoming 
and $m$ outgoing boundaries, and with  morphisms
 the appropriate mapping class groups  
(see Fig.  \ref{catS}). The category $\ES$ is defined in such a  
way that  $\ES(n,1) = \M(n)$. 
So an object of $M(*)$ is a morphism from 0 to 1 in $\ES$, and  
thus defines  
a 1-simplex in $B\ES$. Hence there is  a natural map 
$$\phi:M(*)\x M(*) \rar \Om B\ES$$
as two elements of $\ES(0,1)$ define a loop in $B\ES$. 
We use Barratt and Eccles' method to give the spectra of deloops explicitly,  
obtaining two sequences of simplicial spaces with 
space of $p$-simplices  
$E^i_p=\G\Ga(S^i\w M^p((M(*)\x M(*))_+))$ and  
$F^i_p=\G\Ga(S^{i-1}\w \Ga^p(B\ES))$, where $\Ga$ is  
the $E_\infty$-operad with $\Ga(k)=E\Si_k$, $\G$ is the group 
completion, and $M^p$ (resp.\ $\Ga^p$) means the functor iterated $p$  
times. There is a map of operads $\M\to\Ga$. 
 
\begin{thm} \label{main}
The adjoint of the map $\phi:M(*)\x M(*) \rar \Om B\ES$  
and the operad map $\M\to \Ga$ induce maps 
$$f_p^i:E^i_p=\G\Ga(S^i\w M^p((M(*)\x M(*))_+)) \rar  
        F^i_p=\G\Ga(S^{i-1}\w \Ga^p(B\ES))$$  
for $i\ge 1$ and $p\ge 0$, 
which can be rectified into an equivalence of spectra 
$$(f')^i: (E')^i \sta{\simeq}{\rar} (F')^i, $$  
where $E'$ and $F'$ are spectra equivalent to $E$ and $F$ respectively.  
\end{thm} 
The maps $f^i_p$ are almost simplicial maps in the sense that they  
satisfy all the simplicial identities except for $\de_pf^i_p$  
which is only homotopic to $f^i_{p-1}\de_p$. 
The map $f'$ is a rectification of $f$ in the sense that the 
equivalence $E'\simeq E$ and $F'\simeq F$ is natural with 
respect to $f$ and $f'$.

\vs

In Section \ref{diagsec}, we describe the method of rectification of 
diagrams which will be used in the proof.
In Section \ref{opmcg}, we give the construction of the 
 operad $\M$, following \cite{T},  
spelling out the details needed further on in the text 
and correcting a minor mistake. 
We also describe  
the spectrum of deloops of $\Z\x B\Ga_\infty^+$ produced by $\M$. 
In Section \ref{disjoint}, we give a description of the 
 category $\ES$, adapted to our needs, and    
 produce an actual map inducing the  
equivalence $\Om B\ES\simeq\Z\x B\Ga_\infty^+$. 
In the appendix, we related this map to Tillmann's original proof of this 
equivalence. 
In Section \ref{comp}, we compare the two infinite loop space structures: 
we  construct  the map, rectify it, and show that it induces an  
equivalence of spectra.

\vs

We will be working most of the time in $\Top_*$, the category of  
pointed topological spaces. Most of our spaces are realization of  pointed 
simplicial spaces. We will use the notation $X_\bullet$ for a simplicial  
space and $X=|X_\bullet|$ for its realization.

\vs

I would like to thank my 
supervisor Ulrike Tillmann for her support and encouragements, 
and Bill Dwyer for a very illuminating discussion  
about higher homotopies.

\section{Rectification of diagrams}\label{diagsec}

\begin{figure}[ht]
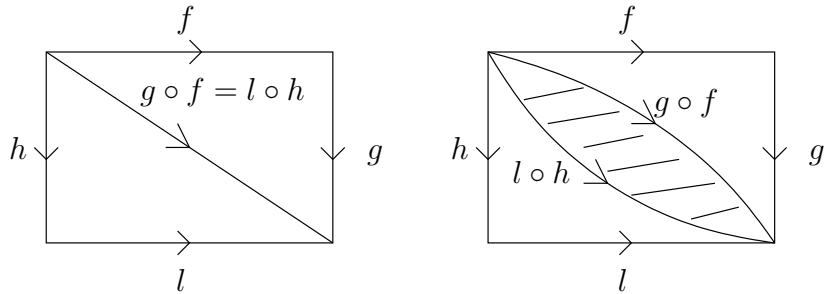

\begin{center}
\input diagram.pstex_t
\caption{Commutative and homotopy commutative square}
  \label{diagram}
\end{center}
\end{figure}

Suppose we have a diagram of spaces and maps (of whatever shape, possibly 
infinite) which 
commutes only up to homotopy. If there are  
higher homotopies, it is possible to rectify it    
to a strictly commutative diagram,  
equivalent to the one we started with (in a sense to be made precise). 
We give here a method which is a special 
case of a theory treated by  Dwyer and Kan in \cite{DK}.
 The same construction was used by Segal in \cite{Seg}.
The idea is to look at a commutative diagram as a functor from a
discrete category $\D$ to $\Top_{(*)}$, the category of (pointed) 
topological space,   and a homotopy 
commutative diagram as a functor from a category $\tD$ to $\Top_{(*)}$,
where the spaces of maps in $\tD$ are ``thicker'' than in $\D$ 
(see Fig.  \ref{diagram} for the case of a square).
As long as the morphism spaces in $\tD$ are homotopy equivalent to 
the corresponding ones in $\D$, a rectification can be constructed.

We give here a precise description of the rectification in the unpointed 
case and prove a strong naturality statement 
(Proposition \ref{equip}), as we need to know more about the 
rectification than what
can be found in \cite{DK} or \cite{Seg}.
The pointed case is done similarly.

\vs

Let $\D$ be a discrete category and let $\tD$ be a category
enriched over $\Top$ with  the same objects as $\D$ and 
such that there is a functor (path components functor) 
 $$p:\tD\to\D,$$ 
which is the identity on objects and induces a homotopy equivalence
$$\tD(x,y)\simeq \D(x,y)$$ 
for each pair of objects $x,y$. So $\tD$ has a contractible space of 
morphisms over each morphism in $\D$ and $p$ is the projection.
There is an induced functor 
\[\Top^\D \sta{p^*}{\rar} \Top^{\tD},\]
from {\em $\D$-diagrams} to {\em $\tD$-diagrams}. 
There is also a functor in the other direction:
\[ \Top^\D \sta{p_*}{\lar} \Top^{\tD},\] 
where $p_* F$ is defined on an object $x$ of $\D$  
as the realization of a simplicial space 
whose $n$th space is 
\[(p_* F)(x)_n= \coprod_{y_0,\dots,y_n\in\ Ob\tD}
          F(y_0)\x\tD(y_0,y_1)\x\dots\x\tD(y_{n-1},y_n)
          \x\D(p(y_n),x).\]

Two functors $F$ and $G$ are said to be {\em equivalent}, denoted 
 $F\simeq G$, if there is a zig-zag
of natural transformations $F\lar F_1\rar \dots\lar F_k\rar G$ which 
induces homotopy equivalences on objects.
\begin{prop}\label{equip}
There is an equivalence of functors 
\[p^* p_* F \simeq F\]
for any $F$ in $\Top^{\tD}$, which is natural in $F$. 
\end{prop}

As $\D$ and $\tD$ have the same objects, this  
means in particular that $p_*F(x)\simeq F(x)$ for any object $x$. 
The functor $p_*F$ is the {\em rectification} of $F$.

\begin{proof}
To prove the Proposition, we will give an explicit sequence 
of natural transformations giving the equivalence and show that they 
are moreover natural with respect to $F$.

For a functor $F:\tD \to \Top$, 
define the functor $\overline{p^*p_*F}$ from $\tD$ to $\Top\,$ simplicially by 
\[(\overline{p^*p_* F})(y)_n= \coprod_{y_0,\dots,y_n\in \ Ob\tD}
          F(y_0)\x\tD(y_0,y_1)\x\dots\x\tD(y_{n-1},y_n)
          \x\tD(y_n,y).\] 
Then there are natural transformations
\[p^*p_*F \lar \overline{p^*p_*F} \rar  F\]
inducing equivalences 
$p^*p_*F(y) \simeq \overline{p^*p_*F}(y) \simeq F(y)$,
for all $y$ in $\tD$.
The natural transformation $\overline{p^*p_*F}\to p^*p_*F$ is
induced by the projection functor $\tD\to \D$ which is a homotopy
equivalence on the space of morphisms, so it clearly induces an 
equivalence.
Now $\overline{p^*p_*F}$ is of the form of a two-sided bar
construction $B(F,D,DX)$.  
May gives an explicit simplicial homotopy for the equivalence 
$B(F,D,DX)\simeq FX$ (\cite{G}, Proposition 9.9). 
It can be adapted to our case.
Indeed, consider the inclusion 
$i:F(y)\hookrightarrow \overline{p^*p_*F}(y)$ defined by 
$a \mapsto (a,id_y)\in F(y)\x\tD(y,y)$, and the 
obvious evaluation map 
$d:\overline{p^*p_*F}(y)\to F(y)$. 
Clearly, $d\circ i=id$. We are left to show that $i\circ d$ is
homotopic to the identity. 
The simplicial homotopy is given explicitly on $q$-simplices by 
\[h_i=s_q\dots s_{i+1}\circ\eta\circ\de_{i+1}\dots\de_q,\] 
for $i=0,\dots,q$, 
where $\eta:(\overline{p^*p_*F})_{i}\to(\overline{p^*p_*F})_{i+1}$ is
defined by adding $id_x\in\tD(x,x)$ on the right of the simplex.

\vs

The first natural transformation is clearly natural in $F$.
The second is natural in $F$ as the diagram 
\[\xymatrix{F(y_0)\times\tD(y_0,y)\ar[d] \ar[r] & F(y) \ar[d]\\
F'(y_0)\times\tD(y_0,y) \ar[r] & F'(y)
}\] 
commutes because a map of functors is itself a natural transformation.
\end{proof}
Note that the inclusion $i:F(y)\hookrightarrow \overline{p^*p_*F}(y)$ 
is not natural in $y$.

\section{The mapping class groups operad}\label{opmcg}

In this section, we describe Tillmann's operad $\M$, correcting a  
minor mistake from \cite{T} in the construction.  
We also give explicitly 
the spectrum of deloops of $\Z\x B\Ga_\infty^+$ produced by $\M$.  
 
\vs

Let $F_{g,n+1}$ denote an oriented surface of genus $g$ with $n+1$ 
boundary components. One of the boundary components is marked; we 
call the $n$ other components {\em free}. Each free boundary component
$\del_i$ comes equipped with a collar, a map from $[0,\epsilon)\x S^1$
to a neighborhood of $\del_i$; for the marked boundary component, there 
is a map from $(\epsilon,0]\x S^1$ to a neighborhood of the boundary. 
Let $\Dif(F_{g,n+1};\del)$ be the group of orientation preserving 
diffeomorphisms which fix the collars, and let 
\[\Ga_{g,n+1}=\pi_0(\Dif(F_{g,n+1};\del))\] 
be its group of components, the associated {\em mapping class group}.

We want to construct a topological operad $\M$ with space of $k$-ary 
operations
$$\M(k)\simeq\coprod_{g\ge 0}B\Ga_{g,k+1}$$
and composition maps induced by gluing surfaces.  
To make gluing associative, one has to  replace the groups $\Ga_{g,k+1}$ 
by equivalent groupoids.

\subsection{Construction of the operad}\label{egn1}

\begin{figure}[ht]
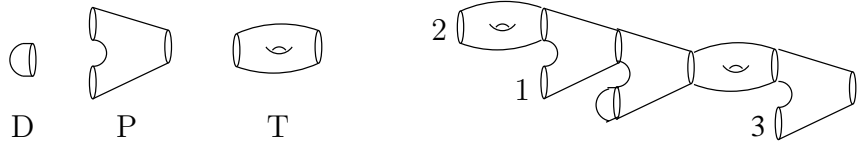

\begin{center}
\input basics.pstex_t
\caption{Building blocks of $\e_{g,k,1}$ and element of $\e_{2,3,1}$}
  \label{basic}
\end{center}
\end{figure}

Pick  a disc $D=F_{0,1}$, a pair of pants surfaces $P=F_{0,3}$ 
and a torus $T=F_{1,2}$ with two boundary components, all with fixed  
collars of the boundary components (see Fig.  \ref{basic}).
Define a groupoid $\e_{g,n,1}$ with objects $(F,\s)$, where 
$F$ is a surface of type $F_{g,n+1}$ constructed from $D,P$ and $T$ 
by gluing the marked boundary of one surface to one of the free 
boundaries of another using the given parametrization, and $\s$ is 
an ordering of the $n$ free boundary components (see Fig.  \ref{basic}). 
Note that each boundary component of $F$ comes equipped with a collar.
The  morphisms from $(F,\s)$ to $(F',\s')$ are the homotopy 
classes $\Ga(F,F')=\pi_0 \Dif(F,F';\del)$ of orientation 
preserving diffeomorphisms preserving the collars and the ordering 
of the boundaries.
The group $\Si_n$ acts freely on $\e_{g,n,1}$  
by permutting the labels and  $B\e_{g,n,1}\simeq B\Ga_{g,n+1}$.

\vs

Gluing of surfaces induces now an associative 
operation on the categories. Hence we  have maps 
on the classifying spaces 
\[\ga:B\e_{g,k,1}\x B\e_{h_1,n_1,1}\x\dots\x B\e_{h_k,n_k,1}\to
      B\e_{g+h_1+\dots+h_k,n_1+\dots+n_k,1}\]
induced by gluing the $k$ last surfaces to the first one according to the  
labels of its free boundaries.  
These maps are associative and $\Si$-equivariant.  
However, $\{\coprod_{g\ge 0}B\e_{g,n,1}\}_{n\in\N}$ does not
precisely form an operad yet as there is no unit.
We will apply a quotient construction on 
the categories $\e_{g,n,1}$ which will both provide a unit and   
make the product induced by the pair of pants associative and 
unital.

\subsubsection{Quotient construction}
\label{quotient}
  
To make the multiplication induced by the pair of pants associative,  
we need to identify subsurfaces of the form $\ga(P;\_,P)$ to  
subsurfaces of the form $\ga(P;P,\_)$. For the unit, we need to  
identify $\ga(P;\_,D)$ and $\ga(P;D,\_)$ to a circle. This circle will  
also be a unit for the operad. 
In \cite{T}, Tillmann does a quotient construction by picking  
morphisms $\phi_1:\ga(P;D,\_)\to\ga(P;\_,D)$ and  
$\phi_2:\ga(P;\_,P)\to\ga(P;P,\_)$ and uses composition of these morphisms  
to identify the surfaces. She then chooses an identification  
of $\ga(P;\_,D)$ to the circle and repeats the process. This is not  
precisely correct as any choice of $\phi_1,\phi_2$ would not yield  
associative operad maps on the quotient categories. It will work only  
with the canonical choice which is the ``identity".  
We prove here that this canonical choice exists and that it  
makes the quotient  
construction possible. We do both quotient constructions at once. 
 
\vs 

Define $\ES_{0,1,1}$ to be the category with one object (thought of as the 
circle) and with $\Z$ as set of morphisms (thought of as the Dehn twists 
around that circle, which will make sense when we glue the 
circle to another surface). 
For $(g,n)\neq (0,1)$, define $\ES_{g,n,1}$ to be the 
full subcategory of $\e_{g,n,1}$ with set of objects all surfaces which do 
not contain subsurfaces of the form $\ga(P,\_,P)$, $\ga(P,\_,D)$, and 
$\ga(P,D,\_)$.

\begin{claim} For each object $F$ in $\e_{g,n,1}$, there is a unique object 
$\overline{F}$ in $\ES_{g,n,1}$  
 obtained from $F$ by a sequence of the following 
operations: replacing a subsurface $\ga(P,\_,P)$ by $\ga(P,P,\_)$, and 
collapsing a subsurface of the form $\ga(P,\_,D)$ or $\ga(P,D,\_)$ to 
a circle. 
In particular, the gluing operation on objects of 
$\coprod\ES_{g,n,1}$ defined 
by   $F\sq_\ES G:=\overline{F\sq G}$, using the gluing $\sq$ defined  
on $\e_{g,n,1}$, is associative. 
\end{claim}

\begin{claim}\label{claim2}
For each $F$ in $\e_{g,n,1}$, one can define a morphism 
      $$\phi_F: F\to \overline{F}$$ 
in $\e_{g,n,1}$, such that all diagrams of the form  
\[\xymatrix{
 & F\sq \overline{G\sq H} \ar[rr]^{\phi_{F\sq \overline{G\sq H}}} 
 && \overline{F\sq \overline{G\sq H}} \ar[dd]^= \\
F\sq G\sq H \ar[ur]^{{\rm id}\sq\phi_{G\sq H}} 
            \ar[dr]_{\phi_{F\sq G}\sq {\rm id}}& & \\
 & \overline{F\sq G}\sq H \ar[rr]_{\phi_{\overline{F\sq G}\sq H}}
 && \overline{\overline{F\sq G}\sq H}
}\]
commute.
\end{claim}

\noindent
{\it Proof of Claim 1.} 
Consider first an element $(F,\s)$ of $\e_{0,n,1}$, i.e. a surface $F$ 
of genus 0 together with a labeling $\s$ of the free boundary 
components ($\s$ is a permutation of the canonical labeling). 
The operations allowed do not change the 
number of free boundaries, nor does it permute the boundaries. In 
$\ES_{0,n,1}$ there is only one surface with labeling $\s$, and this 
surface can clearly be obtained from $F$ by a finite sequence of  
the prescribed moves. 
This surface is $\overline{F}$. 

For a general surface $F$ with labeling $\s$, the moves only affect 
the subsurfaces of $F$ built out of $P$'s and $D$'s. Each of the 
maximal such subsurface has a unique image in the relevant 
$\ES_{0,n_i,1}$. So $\overline{F}$ is the unique surface obtained by 
transforming each of those subsurfaces of $F$.

Finally, gluing as defined in the Claim is associative by the 
uniqueness of the representative of $F\sq G\sq H$ in $\ES_{g,n,1}$.
\cqfd

\vs

\noindent
{\it Proof of Claim 2.} 
Fix three non-intersecting curves on the pair of pants $P$, 
from 0 of the marked boundary to 
0 of the first free boundary, from $\pi$ of the this boundary to 0 of the 
second free boundary, and from $\pi$ of this boundary to $\pi$ of 
the marked boundary (where we think of $S^1$ as parametrized by 
$[0,2\pi[$). This divides the pair of pants into two discs. 
Fix also a curve on the disc $D$, from 0 to $\pi$ of its boundary. 

Now any surface built out of $P$'s and $D$'s comes equipped with a 
system of curves dividing the surface into two discs. These curves run from 
0 of the marked boundary to 0 of the first (in the canonical ordering) 
free boundary, then from $\pi$ of that boundary to 0 of the next, and so on 
until one goes back to $\pi$ of the marked boundary. 
Choose a map $F\to\overline{F}$ which sends the curves of $F$ to the  
corresponding ones in $\overline{F}$. As the curves divide the surfaces  
into discs, 
 by the ``Alexander trick''  
this map is unique up to isotopy. It is, up to isotopy,  
 the identity on the discs.  
Now define 
$\phi_F$ to be the component of this map in $\Dif(F,\overline{F})$.  

For a general surface $F$, define $\phi_F$ to be  
 the map defined by the above on each maximal subsurface of 
$F$ built out of $P$'s and $D$'s, and the identity on the tori. 
The diagram in the Claim commutes by the Alexander trick. 
\cqfd

\vs

Now one can define an operad structure on the categories 
$\ES_{g,n,1}$.
The new structure maps $\overline{\ga}$ are defined on objects by 
taking the unique representative of the image of $\ga$ on the 
surfaces:
\[\overline{\ga}(F,G_1,\dots,G_k)=\overline{\ga(F,G_1,\dots,G_k)}.\]
On morphisms, 
\[\overline{\ga}(f,g_1,\dots,g_k)=
               \phi_{H'}\ga(f,g_1,\dots,g_k)\phi_H^{-1},\] 
where $H$ and $H'$ 
are the images by $\ga$ of the sources and target surfaces of the 
maps $f$ and $g_i$. 
The associativity of 
$\overline{\ga}$ on objects follows from the associativity of 
gluing. On morphisms, it follows from the commutativity of the 
diagram in Claim \ref{claim2}.

Define the operad $\M$ by  
$$\M(k)=\coprod_{g\ge 0}B\ES_{g,k,1}$$ 
with structure maps induced by $\overline{\ga}$. 
Note that there is a map of operads
\[\M \sta{\pi}{\rar}  \Ga, \]
where $\Ga$ is the infinite loop space operad 
with $k$th space $\Ga(k)=E\Si_k$, the classifying space of the 
translation category of the symmetric groups \cite{BE}, and 
$\pi$ is given on objects by the projection to the
labels.

\subsection{Infinite loop space structure}\label{infiniteM}

Let $\mathcal{G}$ denote the group completion functor from the 
category of 
monoids to the category of 
groups 
and let $(M,\mu_M,\eta_M)$ 
denote the monad associated to the operad $\M$ \cite{G}. 
The pair of pants multiplication induces a monoid structure on  
$\M$-algebras.  
$M(*)$ is an $\M$-algebra and 
$$\G M(*)\simeq \G(\coprod_{g\ge 0}B\Ga_{g,1})
                                   \simeq \Z\x B\Ga^+_\infty.$$

Tillmann showed that 
if $X$ is an $\M$-algebra, $\G(X)$ is weakly homotopy equivalent to 
an infinite loop space.
In particular, $\G M(*)\simeq \Z\x B\Ga^+_\infty$ is an infinite 
loop space. 
The $i$th deloop of the group completion of an $\M$-algebra  
$X$ obtained in \cite{T} are defined as the simplicial space  
with $\G\F(S^i\w M^p(M(*)\x X))$ as space of $p$-simplices, where 
$\G\F(Y)$ is the fibre of the map $\G M(Y)\to \G M(*)$.  
As $\G\F(X)\simeq\G\Ga(X)$, one can show that it is equivalent to 
the simplicial space   with $p$-simplices 
\[E^i_p=\G\Ga(S^i\w M^p((M(*)\x M(*))_+)\] 
and that this equivalence induces an equivalence of spectra.  
We work with an added  basepoint 
for technical reasons when constructing the map (this does not bring any  
new component because of the group completion). 

We describe next the simplicial structure on $E^i_\bu$.   
Let $\mu_M,\eta_M$ and $\mu_\Ga,\eta_\Ga$ denote the product and  
unit maps of the monads $M$ and $\Ga$  and 
let $\te$ denote the $M$-algebra structure map of $(M(*)\x M(*))_+$.  
Let $\pi:\M\to \Ga$ 
be the projection of operads. 
 
For any operad $\p$, there is an assembly map $a:A\x P(X) \to P(A\x X)$,  
sending an element $(a,p,x_1,\dots,x_k)$ to $(p,(a,x_1),\dots,(a,x_k))$.  
As $\Ga(*)=\{*\}$, the sequence of maps   
\[\Ga(A\x M(X))\sta{a}{\rar} \Ga(M(A\x X))  
\sta{\pi}{\rar} \Ga\Ga(A\x X) \sta{\mu_\Ga}{\rar} \Ga(A\x X)\]  
induces a map on smash products 
\[\la:\Ga(A\w M(X))\rar \Ga(A\w X).\]  
The simplicial structure on $E^i_\bu$, $i\ge 0$, is defined as follows: 
\[\begin{tabular}{ll} 
$\de_0=\G(\la): E_p^i \to E_{p-1}^i$; & \\ 
$\de_i=\G\Ga(S^i\w M^{i-1}(\mu_M)): E_p^i \to E_{p-1}^i$ & 
                                                 for $1\le i<p$;\\ 
$\de_p=\G\Ga(S^i\w M^{p-1}(\te)): E_p^i \to E_{p-1}^i$; & \\ 
$s_i=\G\Ga(S^i\w M^{i}(\eta_M): E_p^i \to E_{p+1}^i$ &   
                                                for $0\le i\le p$. 
\end{tabular}\]

Let $\uF=F_0\sta{f_1}{\lar}\dots\sta{f_q}{\lar}F_q$ and 
$\uG=G_0\sta{g_1}{\lar}\dots\sta{g_q}{\lar}G_q$ denote  
 elements of $M(*)$. Let $D$ denote the 0-simplex 
represented by the disc.

\begin{prop}\label{equiM} 
Let $\phi:M(*)\to |\G\Ga(M^\bu((M(*)\x M(*))_+))|$ be the map which 
sends $\uF$ to  
the 0-simplex $(1,D,\uF)\in\G\Ga((M(*)\x M(*))_+)$.  Then there is  
a commutative diagram  
$$\xymatrix{M(*) \ar[rr]^{\phi\ \ \ \ \ \ \ \ \ \ \ \ } \ar[d]  
                     & & |\G\Ga(M^\bu((M(*)\x M(*))_+))|\\ 
\G M(*) \ar[urr]_\simeq &  & 
}$$ 
\end{prop} 
 
This can be proved by studying the map of fibrations which Tillmann  
uses to prove the equivalence $\G X\simeq \G\F(M^\bu(M(*)\x X))$.

\section{The cobordism category}\label{disjoint}
 
In this section, we first set up a variant of Tillmann's cobordism  
category $\ES$. In \cite{T2}, the morphism spaces are categories  
similar to the categories $\e_{g,n,1}$ defined in Section \ref{egn1}.    
Our version of $\ES$ is obtained by applying   
the quotient construction of Section \ref{quotient} to Tillmann's  
$\ES$. This version is more amenable to the comparison of the  
two infinite loop space structures.  
We then construct an explicit equivalence 
$\Z\x B\Ga_\infty^+\sta{\simeq}{\rar}\Om B\ES$. We show in the appendix 
how this proof relates to Tillmann's proof.

\vs

The objects of $\ES$ are the natural numbers 
$0,1,2,\dots$. The morphism space $\ES(n,m)=B\ES_{g,n,m}$, where  
$\ES_{g,n,m}$ is a category whose objects are surfaces built out of  
$P$, $T$ and $D$ as in the case of the operad $\M$ but allowing  
disjoint union of surfaces (with a component-wise quotient construction) 
and labeling both the inputs and the outputs 
(see Fig.  \ref{catS}).
The morphisms of $\ES_{g,n,m}$ are homotopy classes of diffeomorphisms
preserving the orientation, the collars and the ordering of the 
boundary components. 
 In particular $\ES(k,1)=\M(k)$. 
Also, a morphism from $n$ to $m$ has exactly $m$ components. 
The only morphism to 0 is the identity in $\ES(0,0)$.  
\begin{figure}[ht]
\begin{center}
\input catSs.pstex_t
\caption{Morphism from 5 to 2 in $\ES$}
  \label{catS}
\end{center}
\end{figure}
Note that $\ES(n,n)$ contains the 
symmetric group, represented by disjoint copies of the circle with 
labels ``on each side''. 
Composition in $\ES$ is induced by gluing the surfaces according to 
the labels, which can be done using the structure maps
of the operad $\M$ on each component. 
Disjoint union of surfaces induces a symmetric monoidal 
structure on $\ES$. As $B\ES$ is connected, it is an infinite loop space.

\vs

We use Barratt and Eccles' machinery \cite{BE} to produce 
 the deloops of $\Om B\ES\simeq\Z\x B\Ga^+_\infty$.
The space of $p$-simplices of the $i$th deloops is given by
\[F^i_p=\G\Ga(S^{i-1}\w \Ga^p(B\ES))\ \ \textrm{for}\ i\ge 1.\]
The simplicial structure of $F^i_\bu=\G\Ga(S^{i-1}\w \Ga^\bu(B\ES))$
is  similar to the one of
$E_\bu^i$, which is given in detail in Section \ref{infiniteM}.

Note that $\ES(0,1)=\M(0)=M(*)$. Recall that 
the pair of pants multiplication 
induces a monoid structure on this space and that
$\G M(*)=\G(\ES(0,1))$ denotes its group completion, which is  
homotopy equivalent to $\Z\x B\Ga^+_\infty$.
In the following Proposition, we use the fact that a morphism in $\ES$  
is a 1-simplex in $B\ES$, and hence two morphisms from 0 to 1 define a loop  
in $B\ES$. 

\begin{prop}\label{equiS}
Define  $\psi:\ES(0,1) \to \Om B\ES$
by  $\psi(\uF)=\uF\Big(^1_0\Big)D$ is the 
loop from 0 to 1 along the morphism defined by $\uF$ followed by the 
morphism defined by the disc $D$ taken backwards.
Then there is a homotopy commutative diagram 
\[\xymatrix{M(*)=\ES(0,1) \ar[rr]^\psi \ar[d] & & \Om B\ES\\
\G M(*)=\G\ES(0,1). \ar[urr]_\simeq & &
}\]
\end{prop}

\begin{proof}
Consider the diagram 
$$\xymatrix{
\ES(0,1) \ar[rrr]^\psi \ar[d]\ar[dr]^f & 
        & 
        & \Om B\ES \ar[d]^=\\
\G\ES(0,1) \ar[r]^\simeq & \Om B(\ES(0,1))\ar[r]_g^\simeq 
   & \Om_{1,1}(B\ES) \ar[r]_h^\simeq  & \Om B\ES
}$$
where $\Om_{1,1}(B\ES)$ denotes the space of paths from 1 to 1 
in $B\ES$. 
The map $f$ is Quillen's group completion map, which sends  
an element $\uF$ of the monoid $\ES(0,1)$ to the loop it defines in its 
classifying space (as a monoid). 
To define $g$, consider 
the map $\ES(0,1)\to\ES(1,1)$ induced by gluing a pair 
of pants: $\uF \mapsto \uF\sq P$ (glue the surface to the left leg, 
i.e.  compose $\uF\coprod S^1$ with $P$ in $\ES$). 
This induces a functor from the monoid $\ES(0,1)$ to the category $\ES$. 
Indeed, one can check that  
the pairs of pants multiplication is mapped to composition 
in $\ES$. Tillmann showed in \cite{T1} (Proposition 4.1) that $g$ induces a 
homotopy  equivalence 
on the classifying spaces. 
Finally, the map $h$ is defined by precomposing with the path from 0 to 1 
(1-simplex) given by the disc and postcomposing by the same path 
taken backwards. 

The diagram commutes up to homotopy. Indeed, starting with $\uF$ in 
$\ES(0,1)$, the loop obtained by following the bottom of the diagram 
is going along $D$ from 0 to 1, then $\uF\sq P$ from 1 to 1 and lastly 
$D$ again backwards from 1 to 0. This path is homotopy equivalent to 
$\psi(\uF)$ as $D\circ(\uF\sq P)=\uF$ in $\ES$, which means that there 
is a 2-simplex in $B\ES$ providing the required homotopy.
\end{proof}

Tillmann's 
proof that $\Om B\ES \simeq \Z\x B\Ga^+_\infty$ is by showing that 
$\Z\x B\Ga^+_\infty$ is equivalent to a homotopy fibre which is known to 
be $\Om B\ES$. We will show in the appendix that the map described 
above is a natural choice in this context to make the equivalence 
explicit. This leads to a more natural proof of the Proposition.

\section{Comparison of the two structures}\label{comp}

In this section, we compare the two infinite loop space structures.  
We first construct in \ref{mapfp} the maps $f^i_p:E^i_p\to F^i_p$  
of Theorem \ref{main}.   
In \ref{rectfp}, we rectify these maps to simplicial maps  
between simplicial spaces $(E')^i_\bu$ and $(F')^i_\bu$  
equivalent to $E^i_\bu$ and $F^i_\bu$. 
In \ref{specfp}, we show that this rectification provides 
 a map of spectra which  in  
\ref{equfp} is shown to be an equivalence.   
Theorem \ref{Dthm}, Theorem \ref{spectra} and Theorem \ref{final} combine to prove  
the main Theorem \ref{main}.

\subsection{Construction of a map}\label{mapfp}

We want to construct a map from
$E_p^i=\G\Ga(S^{i}\w M^p((M(*)\x M(*))_+))$, the space of $p$-simplices 
of the $i$th deloop of $\G M(*)$,  to 
$F_p^i=\G\Ga(S^{i-1}\w \Ga^p(B\ES))$, the $p$-simplices of the 
$(i\!-\!1)$st deloop of $B\ES$, for $i\ge 1$. 
By construction, an element of $M(*)$ 
is a morphism from 0 to 1 in the category $\ES$, and hence a  
1-simplex in its classifying space $B\ES$. 
In particular, two such elements define a loop in $B\ES$.  
\begin{figure}[ht]
\begin{center}
\input map.pstex_t
\caption{Map from $S^1\w (M(*)\x M(*))_+$ to $B\ES$}
  \label{map}
\end{center}
\end{figure}
The map obtained using this remark is naturally bisimplicial: 
\[\tilde{\phi}_{p,q}: S^1_p \x(M_q(*)\x M_q(*)) \rar B_{p,q}\ES,\]
where $S^1$ is viewed as a simplicial space with two 0-simplices 
$x_0,\, x_1$   
and two non-degenerate 1-simplices $y_1,\, y_2$ (see Fig.  \ref{map}),  
and $B\ES$ is viewed as a bisimplicial space  
with the second simplicial dimension coming  
from the simplicial structure of its morphism spaces. \\
\begin{tabular}{llll}
Define & $\tilde{\phi}_{0,q}(x_i,\uF,\uG)=$ & $i$  & for  $i=0,1$;\\
 & $\tilde{\phi}_{1,q}(y_1,\uF,\uG)=$ &  $0 \sta{\uF}{\rar} 1$; & \\
 & $\tilde{\phi}_{1,q}(y_2,\uF,\uG)=$  & $0 \sta{\uG}{\rar} 1$. &
\end{tabular}\\
This induces a map 
$\phi: S^1 \w(M(*)\x M(*))_+ \rar B\ES$.
The base point of $M(*)\x M(*)$ is $(D,D)$ and its image under 
$\tilde{\phi}$ is a contractible loop, but it is not actually the 
trivial loop. This is why we work with $(M(*)\x M(*))_+$ rather than 
$M(*)\x M(*)$.

Recall from \cite{BE} 
that for the monad $\Ga$ there is an assembly map 
\[a:A\w \Ga(X) \to \Ga(A\w X),\]
defined by 
$a(y,\underline{\s},x_1,\dots,x_n)=(\underline{\s},[y,x_1],\dots,[y,x_n])$. 
This also induces a map $\A\w\G\Ga(X)\to \G\Ga(A\w X)$. 
Combining $a$, $\phi$ and the operad map $\pi:\M\to\Ga$, we get a map
\[\xymatrix{
\G\Ga(S^{i}\w M^p((M(*)\x M(*))_+)) \ar[rr]^{\ \ \ f^i_p} 
  \ar[d]_{\G\Ga(S^i\w \pi^p)} & &  \G\Ga(S^{i-1}\w \Ga^p(B\ES))\\
\G\Ga(S^{i}\w \Ga^p((M(*)\x M(*))_+)) \ar[d]_{\G\Ga(S^{i-1}\w a)} 
   & & \\ 
\G\Ga(S^{i-1}\w \Ga^p(S^1\w (M(*)\x M(*))_+)). 
       \ar[uurr]_{\ \ \ \G\Ga(S^{i-1}\w\Ga^p(\phi))}  & &
}\]

\begin{prop}\label{almost}
Let $\de_j$ and $s_j$ denote the boundary and degeneracy maps of both 
the simplicial spaces  $E^i_\bu$ and $F^i_\bu$. Then the maps 
$f_p^i:E^i_p \to F^i_p$ for $i\ge 1$ and $p\ge 0$ satisfy 
\[\de_j f^i_p = f^i_{p-1}\de_j \hspace{1cm} \textrm{for }\ 0\le j<p\]
and
\[s_j f^i_p = f^i_{p+1}s_j \hspace{1cm} \textrm{for }\ 0\le j\le p.\]
\end{prop}

\begin{proof}
The boundary maps $\de_j$, for $0\le j<p$, and the degeneracies 
for $E^i_\bu$ and $F_\bu^i$ 
are defined in terms of the operad structure maps of $\M$ and $\Ga$. 
The map $f_\bu$ commutes with all of those maps 
because $f$ maps $\M$ to $\Ga$ via an operad map.
\end{proof}

The last commutation relation necessary to have a simplicial map, 
$\de_pf_p=f_{p-1}\de_p$, is satisfied only up 
to homotopy. We show in the next section that it preserves all the  
higher homotopies and hence that it can be rectified into a  
simplicial map.

\subsection{Rectification of the map}\label{rectfp}

\begin{thm}\label{Dthm}
There exist simplicial spaces $(E')^i_\bu$ and $(F')^i_\bu$ 
equivalent to the simplicial spaces $E_\bu^i$ and $F^i_\bu$ and 
a simplicial map $(f')^i_\bu:(E')^i_\bu\to (F')^i_\bu$ such that 
the following diagram commutes:
\[\xymatrix{E^i_p \ar[r]^{f^i_p} & F_p^i \\
(E')^i_p \ar@{<->}[u]^\simeq \ar[r]^{(f')^i_p} 
                & (F')^i_p\ar@{<->}[u]_\simeq}\]
\end{thm}

To prove this theorem, 
we will use the method described in Section \ref{diagsec}.
We first need to construct the categories $\D$ and $\tD$ relevant to our
situation. 

\vs

Let $\dop$ denote the simplicial category: the objects of $\dop$
are the natural numbers and there are maps $\de_i:p\to p-1$ and 
$s_i:p\to p+1$ for each $i=0,\dots,p$, satisfying the simplicial
identities. So  a $\dop$-diagram  is a 
simplicial space. Any morphism in $\dop(p,q)$ can be
expressed uniquely as a sequence 
$s_{j_t}\dots s_{j_1}\de_{i_s}\dots\de_{i_1}$ with 
$0\le i_s< \dots < i_1\le p$ and 
$0\le j_1< \dots < j_t\le q$ and $q-t+s=p$ \cite{M4}.

Let $\D$ be the category whose $\D$-diagrams are
precisely a couple of simplicial spaces with a simplicial map between
them (see Fig.  \ref{catD2}).
\begin{figure}[ht]
\begin{center}
\input catD2s.pstex_t
\caption{Category $\D$}
  \label{catD2}
\end{center}
\end{figure}
So $\D$ has two copies of the natural numbers as set of
objects, denoted $E_p$ and $F_p$ for $p\in \N$.  The full subcategory of
$\D$ containing all the $E_p$'s is isomorphic to $\dop$. So 
$\D(E_p,E_q)=\dop(p,q)$. Similarly $\D(F_p,F_q)=\dop(p,q)$.
Finally, there is a unique map $f_p\in\D(E_p,F_p)$ and it satisfies the
simplicial identities $\de_i f_p=f_{p-1}\de_i$ and 
$s_i f_p=f_{p+1}s_i$ for $i=0,\dots,p$. 
So any morphism in $\D$ from $E_p$ to $F_q$ can
be written uniquely as a sequence 
$s_{j_t}\dots s_{j_1}\de_{i_s}\dots\de_{i_1}f_p$ where the indices
are as above.

We now define the category $\tD$ 
 in such a way that the data given in \ref{mapfp} will
induce a functor from $\tD$ to $\Top_*$.
$\tD$ has the same objects as $\D$ and we will again denote
them by $E_p$ and $F_p$ for $p\in\N$. 
Also, 
$\tD(E_p,E_q)=\D(E_p,E_q)$ and $\tD(F_p,F_q)=\D(F_p,F_q)$.
To describe $\tD(E_p,F_q)$, we first need to define the 
{\em degeneracy degree} $d(g)$ of a morphism $g\in\D(E_p,F_q)$.
If $g=s_{j_t}\dots s_{j_1}\de_{i_s}\dots\de_{i_1}f_p$ in the above  
notation, then $d(g)$
is  the biggest $k$ such that ${i_k}=p-k+1$, and $d(g)=0$ if no such
$k$ exists. In other words, $d$ counts the number of ``bad'' maps,
i.e. last boundary maps, 
occurring in $g$. 
Define
\[\tD(E_p,F_q)=\coprod_{g\in\D(E_p,F_q)} \Delta_{d(g)},\]
where $\Delta_d=\{(t_0,\dots,t_d)\in\RR^{d+1}|t_i\ge 0, \Si t_i=1\}$ 
is the standard $d$-simplex.
Note that for each $g\in\D(E_p,F_q)$, there is an  inclusion 
\[\tilde{g}:\Delta_{d(g)} \rar \tD(E_p,F_q)\]
whose image is the space of morphisms ``sitting over $g$''.

Recall that 
all the simplicial identities between the
$\de_i$'s and $s_j$'s are satisfied in $\tD$. Also, note that 
there is a unique map 
(0-simplex) between $E_p$ and $F_p$ sitting over the simplicial map
$f_p\in\D(E_p,F_p)$. We denote this map again by $f_p$ and we set 
the relation $s_if_p=f_{p+1}s_i$ for $0\le i\le p$ and
$\de_i f_p=f_{p-1} \de_i$ for $0\le i < p$. 
Because the last relation, when $i=p$,  does not hold in $\tD$, 
there are exactly $d(g)+1$ maps formed of compositions of 
$\de_i$'s, $s_j$'s and an $f_k$ projecting down to 
$g$ in $\D$. We define 
those compositions to be the vertices of $\Delta_{d(g)}$.
More precisely, for $0\le k\le d(g)$, define
\[g_k:=
s_{j_t}\dots s_{j_1}\de_{i_s}\dots
                \de_{i_{k+1}}f_{p-k}\de_{i_k}\dots\de_{i_1}
                          :=\tilde{g}(0,\dots,1,\dots,0),\]
where 1 is in the $k$th position counting backwards.

Now composition in $\tD$ is determined by the vertices of the simplices. 
A map in $\tD(E_p,F_q)$ can only be pre-composed by a map 
in $\tD(E_r,E_p)$ or post-composed by a map in $\tD(F_q,F_s)$. 
Using the identities given above, we know 
how those compositions are defined on the vertices of the simplices 
of $\tD(E_p,F_q)$. We then extend the composition simplicially.

\begin{thm}\label{tD}
Let $\tD$ be the category defined above.
For each $i\ge 1$, there is a functor 
\[L_i: \tD \rar \Top_*\]
such that $L_i(E_\bu)=E_\bu^i$, $L_i(F_\bu)=F_\bu^i$,  
where $E_\bu$ and $F_\bu$ denote 
the two subcategories of $\tD$ isomorphic to $\dop$, and $L_i(f_p)=f_p^i$.
\end{thm}

\begin{proof}
As $E_\bu^i$ and $F_\bu^i$ are simplicial spaces, the restriction of $L_i$ to 
each copy of $\dop$ in $\tD$ is a well-defined functor. 
As the map $f_p^i$ satisfies the identities satisfied by $f_p$ with the 
boundary and degeneracy maps, $L_i$ is also well defined on the 
vertices of the 
simplices of the morphism spaces $\tD(E_p,F_q)$. 
We have to show that we can extend the definition of $L_i$ to 
the whole simplices. 

Let $d=d(g)$ and $g_k$ be the $k$th vertex of the $d$-simplex over $g$ 
as described above.
To simplify notations, let  
$$X=(M(*)\x M(*))_+.$$ 
Because only the first $d$ last boundaries appear in $g$, one can 
factorize the map $g_k$ for any $0\le k\le d$ as \\
$g_k: E_p^i=\G\Ga(S^i\w M^p(X))
\sta{\al}{\rar} \G\Ga(S^{i-1}\w\Ga^{p-d}(S^1\w M^d(X)))$\\
.\hfill $\sta{A(h_k)}{\rar} 
     \G\Ga(S^{i-1}\w\Ga^{p-d}(B\ES))=F^i_{p-d}\sta{\be}{\rar} F_q^i$,\\
where $\al$ is a composition of $\pi:M\to\Ga$ and the assembly map,  
 $\be$ is a composition of boundaries and degeneracies 
  and $A$ is the functor 
$\G\Ga(S^{i-1}\w\Ga^{p-d}(\_\_))$. Both $\al$ and $\be$ are 
 independent of $k$. 

We first construct a $d$-simplex of maps 
 $S^1\w M^d(X)\to B\ES$ having the maps $h_k$ as vertices.
For $0\le l< k\le d$, $h_k$ and $h_l$ are given by the two sides of the 
following diagram:
\[\xymatrix{S^1\w M^d(X)\ar[r]^{\de_{i_l}\!\dots\!\de_{i_1}\ \ \ \ }  &
S^1\w M^{d-k}M^{k-l}(X) \ar[r]\ar[d]^{\de_{i_k}\dots\de_{i_{l+1}}}   
   & \Ga^{d-k}\Ga^{k-l}(B\ES) \ar[d]^{\de_{i_k}\dots\de_{i_{l+1}}} & \\
 & M^{d-k}(S^1\w X)\ar[r]  
& \Ga^{d-k}(B\ES) \ar[r]^{\ \ \ \ \de_{i_d}\!\dots\!\de_{i_{k+1}}}   
 & B\ES
}\]
We think of an element of $M^d(X)$ as being divided into   
 $d+1$ levels $d,\dots,0$: 
level $d$ is the surface (with elements of the mapping class group)  
coming from the $M$ the most on the left,
level $d-1$ is composed of $m_{d-1}$ surfaces coming from the 
second $M$, and so on up to level 0 which is composed of $m_0$ 
elements of $X=(M(*)\x M(*))_+$.  In the above diagram, we start by gluing 
the levels $l,\dots,1$ to level 0. 
The maps $\de_{i_k}\dots\de_{i_{l+1}}$ going down
 glue similarly the levels $k,\dots,l+1$, and 
the image of  $h_k$ is  the disjoint union of the surfaces obtained,  
which is  a couple of morphisms $\uF_{k,0},\uG_{k,0}$ 
from 0 to $m_k$ in $\ES$. 
The map $h_l$ is obtained by following the top of the diagram.  
Once the levels $l,\dots,1$ glued 
to level 0, $h_l$ takes the disjoint union, producing 
two morphisms $\uF_{l,0},\uG_{l,0}$ from 0 to $m_l$ in $\ES$. 
Now the surface $\underline{H}_{k,l}$ obtained by gluing 
together the levels $k,\dots,l+1$ gives a morphism from $m_l$ to 
$m_k$ in $\ES$ which is such that 
$\underline{H}_{k,l}\circ\uF_{l,0}=\uF_{k,0}$ and 
$\underline{H}_{k,l}\circ\uG_{l,0}=\uG_{k,0}$. This produces two 
2-simplices in $B\ES$ and hence a homotopy equivalence 
(see Fig.  \ref{edge}). 

\begin{figure}[ht]
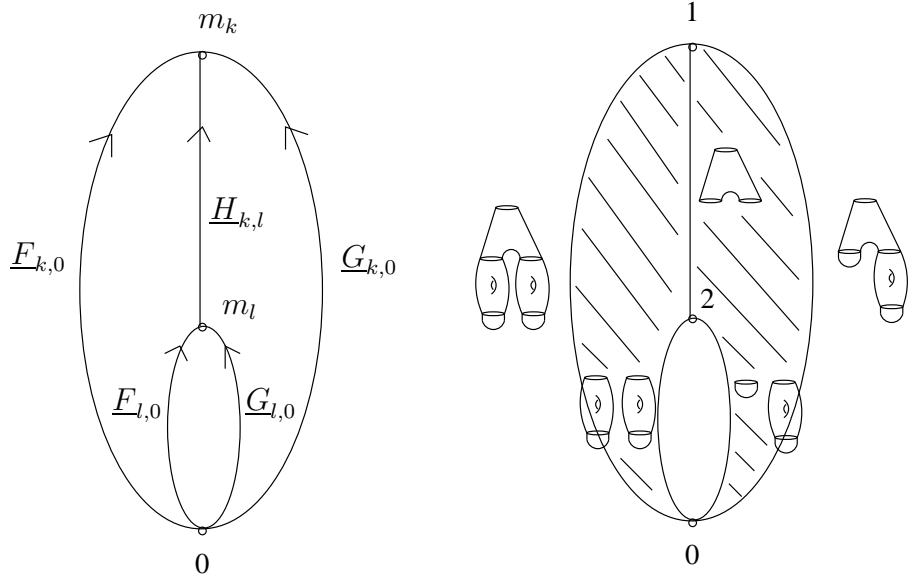

\begin{center}
\input edge.pstex_t
\caption{ Homotopy $H_{k,l}$ }
  \label{edge}
\end{center}
\end{figure}

More precisely, consider  the two $d+1$-simplices of $B\ES$  
\[m_d\sta{\uF_d}{\lar}\dots\sta{\uF_1}{\lar}m_0\sta{\uF_0}{\lar}0
\ \ \ \ \ \ {\rm and}\ \ \ \ \ \  
m_d\sta{\uF_d}{\lar}\dots\sta{\uF_1}{\lar}m_0\sta{\uG_0}{\lar}0,\]
where $\uF_i$ is the disjoint union of the surfaces of level $i$ 
of an element ${\bf F}$ of $M^d(X)$ for 
$i\ge 1$ and 
$(\uF_0,\uG_0)$ is the disjoint union, component-wise, of level 0.  
Then $h_k(S^1,{\bf F})$ is the loop in $B\ES$ from 0 to $m_k$ 
along 
$\uF_{k,0}=\uF_k\circ\dots\circ\uF_0$ and back to 0 along  
 $\uG_{k,0}=\uF_k\circ\dots\circ\uG_0$. The morphism 
$\underline{H}_{k,l}=\uF_k\circ\dots\circ\uF_{l+1}$ induces a 
homotopy between $h_k$ and $h_l$ as explained above.  
Note that  the homotopies $\underline{H}_{k,l}$ 
form the edges of a $d$-simplex.

Define  
\[\tilde{h}:\Delta_d\x S^1\w M^d(X) \rar B\ES\]
for ${\bf F}=(\uF_d,\dots,\uF_1,(\uF_0,\uG_0))$  by 
  
\begin{tabular}{lcl} 
    $\tilde{h}(\Delta_d,x_0,{\bf F})$ & = & 0 \\ 
    $\tilde{h}(\Delta_d,x_1,{\bf F})$ & = &  
$m_d\sta{\uF_d}{\lar}\dots\sta{\uF_1}{\lar}m_0$\\ 
    $\tilde{h}(\Delta_d,y_1,{\bf F})$ & = & 
$m_d\sta{\uF_d}{\lar}\dots\sta{\uF_1}{\lar}m_0\sta{\uF_0}{\lar}0$ \\
$\tilde{h}(\Delta_d,y_2,{\bf F})$ & = &
$m_d\sta{\uF_d}{\lar}\dots\sta{\uF_1}{\lar}m_0\sta{\uG_0}{\lar}0$ 
\end{tabular}\\  
(so $\tilde{h}$ maps $\Delta_d\x I$ to $\Delta_{d+1}$ by collapsing
 $\Delta_d\x\{0\}$). 
This induces  a map 
\[\tilde{g}:\Delta_d\x E^i_p \rar F^i_q\]
by the continuity of the functor $A=\G\Ga(S^{i-1}\w\Ga^{p-d}(\_\_))$,  
which extends the definition already given on the vertices of $\Delta_d$.  
\end{proof}

\noindent
{\em Proof of Theorem \ref{Dthm}.} \hspace{2mm}
The projection $p:\tD\to \D$ induces homotopy equivalences 
$\tD(A,B)\simeq\D(A,B)$. Hence the functor $L_i$ defined in Theorem \ref{tD}
has a rectification 
$L'_i=p_*(L_i): \D\to \Top_*$.  
Denote by $(E')^i_p,\ (F')^i_p$ and $(f')^i_p$ the images of $E_p$, 
$F_p$ and $f_p^i$ via $L'_i$. By definition of the category $\D$, 
 $(E')^i_\bu$ and $(F')^i_\bu$ are simplicial spaces and  
$(f')_\bu^i:(E')^i_\bu\to (F')^i_\bu$ is a simplicial map.

We know that $p^*L'_i\simeq L_i:\tD\to \Top_*$. Now note that 
$p^*(L'_i)\bigr|_{\dop}=(L'_i)\bigr|_{\dop}:\dop\to\Top_*$, 
for each copy of $\dop$ in
$\tD$ and the corresponding one in $\D$, as $p$ is the identity on those 
subcategories. So $L'_i\bigr|_{\dop}\simeq L_i\bigr|_{\dop}$, 
which precisely says that the simplicial
spaces $(E')^i_\bu$ and $(F')^i_\bu$ are equivalent to $E^i_\bu$ and $F_\bu^i$
respectively.  
Lastly, 
as $p^*L'_i(f_p)=(f')^i_p$,
the diagram given in the corollary commutes by naturality of the 
equivalence of functors.
\hfill $\square$

\subsection{Map of spectra}\label{specfp}

The spaces $E^i$ and $F^i$ for $i\ge 1$ form two $\Om$-spectra.  
This means that there are equivalences 
$\epsilon^i:E^i\sta{\simeq}{\rar} \Om E^{i+1}$ and 
$\epsilon^i:F^i\sta{\simeq}{\rar} \Om F^{i+1}$. 
Both adjoints $\overline{\epsilon}^i$ can be expressed simplicially:
\[\epb^i_p: \Si\G\Ga(S^{i}\w M^p((M(*)\x M(*))_+)) \rar
                                 \G\Ga(S^{i+1}\w M^p((M(*)\x M(*))_+))\]
\[\epb^i_p: \Si\G\Ga(S^{i-1}\w \Ga^p(B\ES)) \rar
                                 \G\Ga(S^{i}\w \Ga^p(B\ES)),\]
where both maps are given by the assembly map 
$A\w\G\Ga(X) \to \G\Ga(A\w X)$ described earlier. 
The following Proposition is a direct consequence of the definitions 
of $f^i_p$ and of the spectrum structure maps. 

\begin{prop}\label{spec}
For $p\ge 0$ and $i\ge 1$, the following diagram commutes:\\
$\xymatrix{
\Si\G\Ga(S^{i}\w M^p((M(*)\x M(*))_+)) \ar[rr]^{\Si f^i_p} 
   \ar[d]^{\overline{\epsilon}^i_p}
   & & \Si\G\Ga(S^{i-1}\w \Ga^p(B\ES)) \ar[d]^{\overline{\epsilon}_p^i}\\
\G\Ga(S^{i+1}\w M^p((M(*)\x M(*))_+)) \ar[rr]^{f^{i+1}_p} 
                    & & \G\Ga(S^{i}\w \Ga^p(B\ES)) 
}$
\end{prop}

We want to show that the simplicial spaces $(E')^i$ and $(F')^i$ 
form two spectra equivalent to the original ones and that the 
maps $(f')^i$ form a map of spectra. For this, we need two lemmas.

\begin{lem}\label{nateps}
For each $i\ge 1$, the two sequences of 
maps $\epb^i_p$ induce a natural transformation of functors
\[\epb^i: \Si L_i \rar L_{i+1}.\]
\end{lem}

\begin{proof}
 $\Si L_i$ and $L_{i+1}$ are functors from $\tD$ to $\Top_*$. 
We already know that the $\epb_p^i$'s form a couple of simplicial maps 
and commute with the maps $f_p$ (Proposition \ref{spec}). So we only 
need to check that the $\epb^i_p$ commute with all the homotopies. This 
follows from the fact that the map $\epb^i_p$ is  
 induced by an assembly map $\Si\G\Ga X \to \G\Ga\Si X$, which is natural
in $X$ \cite{BE}. 
\end{proof}

\begin{lem}\label{beta}
For any functor $F:\D\to\Top_*$, there is a natural transformation 
\[\be:\Si(p^*p_*F) \rar p^*p_*(\Si F)\] 
such that the following diagram is commutative: 
$\ \ \ \xymatrix{
\Si(p^*p_*F) \ar[d]_{\be} 
 \ar@{<->}[rr]^{\Si(\simeq)} & & \Si F\\
p^*p_*(\Si F) \ar@{<->}[urr]^{\simeq} & &
}$\\
In particular, $\be$ is an equivalence.
\end{lem}

\begin{proof}
There is a natural map  
$\be_n:(\Si(p_*F)(x))_n \to (p_*(\Si F)(x))_n$  on each simplicial level
as the second space is a quotient of the first. The resulting map $\be$ 
 collapses a contractible subspace of the first simplicial space. 
One then checks that the diagram commutes.
\end{proof}

\begin{thm}\label{spectra}
The spaces $(E')^i$ and $(F')^i$ for $i\ge 1$ form two spectra 
which are equivalent to the spectra $E^i$ and $F^i$, and 
the maps $(f')^i:(E')^i\to(F')^i$ form a map of spectra.
\end{thm}

\begin{proof}
Recall that $p_*$ is a functor $\Top_*^{\tD} \to \Top_*^\D$.
By Lemma \ref{nateps}, we have a natural transformation 
$\Si L_i \sta{\epb}{\rar} L_{i+1}$. Denote by 
$(\Si L_i)'\sta{\epb'}{\rar} L'_{i+1}$ its image under $p_*$.
We need to construct 
maps $\la^i:(E')^i\sta{\simeq}{\rar}\Om(E')^{i+1}$ and 
$\la^i:(F')^i\sta{\simeq}{\rar}\Om(F')^{i+1}$. We define their 
adjoint $\overline{\la}$ simplicially in the following diagram:
\[\xymatrix{
\Si(E')_p^i \ar[rrrrr]^{\Si(f')_p^i} \ar[ddr]^{\be_p} 
            \ar@{<->}[ddrr]^{\Si(\simeq)} 
            \ar@{-->}[ddddr]_{\overline{\la^i_p}} 
 & & & & & \Si(F')_p^i \ar[ddl]_{\be_p} \ar@{<->}[ddll]_{\Si(\simeq)} 
            \ar@{-->}[ddddl]^{\overline{\la^i_p}}\\
 & & & & & \\
& (\Si E')^i_p \ar@{<->}[r]^\simeq \ar[dd]^{\epb'_p} 
   & \Si E^i_p \ar[r]^{\Si f^i_p} \ar[d]^{\epb_p}
   & \Si F^i_p \ar[d]_{\epb_p}
   & (\Si F')^i_p \ar@{<->}[l]_\simeq \ar[dd]_{\epb'_p} & \\
& & E^{i+1}_p \ar[r]^{f_p^{i+1}} & F^{i+1} & & \\
& (E')^{i+1}_p \ar@{<->}[ur]^\simeq \ar[rrr]^{(f')^{i+1}_p}
   & & & (F')^{i+1}_p \ar@{<->}[ul]_\simeq  & 
}\] 
This diagram commutes by Proposition \ref{spec} for the commutation of the 
square in the center, Proposition \ref{equip} 
(naturality in $F$ of the equivalence $p^*p_*F \simeq F$) 
for the commutation of the left and right 
squares, because the equivalence is a natural transformation of
 functors for the top and bottom squares, 
and by Lemma \ref{beta} for the two triangles.
 
Now the adjoint of $\overline{\la}$ are equivalences as 
\[\xymatrix{\Si(E')^i \ar[r]^{\Si(\simeq)} \ar[d]_{\overline{\la}^i} 
                  & \Si E^i \ar[d]^{\epb^i}\\
            (E')^{i+1} \ar[r]^\simeq & E^{i+1}}
\hspace{1cm} \Longrightarrow \hspace{1cm}
\xymatrix{(E')^i \ar[r]^{\simeq} \ar[d]_{\la^i} 
                  & E^i \ar[d]^{\eps^i}\\
           \Om(E')^{i+1} \ar[r]^{\Om(\simeq)} & \Om E^{i+1}}\]
the commutation of the left diagram implies the commutation of the right 
one, and the maps $\eps^i$ are equivalences.
Moreover, this last diagram shows that the equivalences 
$E^i\simeq (E')^i$ and $F^i\simeq (F')^i$ are equivalences of spectra. 

Finally, the commutation of the larger square in the big diagram 
implies that the maps $(f')^i$ form a map of spectra.
\end{proof}

\subsection{Equivalence}\label{equfp} 

So far, we have defined the spectra $(E')^i$ and $(F')^i$ only 
for $i\geqslant 1$, 
i.e. starting with the first deloop of $\Z\x B\Ga_\infty^+$. 
Define $(E')^0:=\Om(E')^1$, $(F')^0:=\Om(F')^1$ and  
\[(f')^0:=\Om (f')^1: (E')^0 \rar (F')^0.\]
Note that $\Om(E')^1\simeq\Z\x B\Ga^+_\infty\simeq\Om(F')^1$.

\begin{thm}\label{final}
The map of spectra $\{(f')^i\}_{i\ge 0}:\{(E')^i\}\rar \{(F')^i\}$ 
is an equivalence.
\end{thm}

\begin{lem}\label{lemeq}
There are maps $\phi_p:M(*)\to \G\Ga(M^p((M(*)\x M(*))_+))$ 
and $\psi:M(*)\to \Om B\ES$ such that the following diagram 
commutes.
\[\xymatrix{
(E')^1_p \ar[rr]^{(f')^1_p} \ar@{<->}[d]_\simeq & & 
    (F')^1_p \ar@{<->}[d]^\simeq \\
E^1_p=\G\Ga(S^1\w M^p((M(*)\x M(*))_+)) \ar[rr]^{f_p^1} 
                    & &  \G\Ga(\Ga^p(B\ES))=F^1_p \ar[d] \\
\Si\G\Ga(M^p((M(*)\x M(*))_+)) \ar[u]  & &  B\ES\\
& \Si M(*) \ar[ul]^{\Si(\phi_p)} \ar[ur]_{\overline{\psi}}  &   
}\]
\end{lem}

\begin{proof}
Let $\uF$ be an element of $M(*)$. 
Define  $\phi_p(\uF):=(1,\dots,1,(\uF,D))$, where $D$ is the disc,  
and let $\psi$ be the map defined in Proposition \ref{equiS} 
which sends $\uF$ 
to the loop in $B\ES$ going from 0 to 1
along $\uF$ and back to 0 along $D$. 
The bottom part of the diagram is easily seen to commute. 
The top part commutes by naturality of the equivalence.
\end{proof}

\noindent
{\em Proof of Theorem \ref{final}.} \hspace{2mm}
The spectra $(E')^i$ and $(F')^i$, for $i\ge 0$ are connective.  
Indeed, they are equivalent to the spectra $E^i$ and $F^i$.  
As the functor $\Ga$ preserves connectedness, 
 $E^i$ and $F^{i+1}$ are connected  
for $i\ge 1$. Moreover, $F^1\simeq B\ES$  is also connected.  
Hence  
it is enough to show that $(f')^0$ is an
equivalence.  
 
Thinking of $M(*)$ as a constant simplicial space, 
Lemma \ref{lemeq} yields a commutative diagram of simplicial spaces 
(with no map from $E^1_\bu$ to $F^1_\bu$). Taking adjoints and using 
Propositions \ref{equiM} and \ref{equiS}
 we get a homotopy commutative diagram
\[\xymatrix{
(E')^0 \ar@{<->}[d]_\simeq \ar[rrrr]^{(f')^0}  & & & &  (F')^0 \\
|\G\Ga(M^\bu((M(*)\x M(*))_+))| 
& &  M(*) \ar[ll]_{\ \ \ \ \ \ \ \phi} \ar[rr]^\psi \ar[d] 
                              & & \Om B\ES \ar@{<->}[u]_\simeq \\
& &  \G M(*) \ar[ull]^\simeq \ar[urr]_\simeq & &
}\] 
(the left triangle commutes only up to homotopy).
Hence $(f')^0$ is a homotopy equivalence.  
\cqfd

\section*{Appendix} {\small

The proof that $\Om B\ES\simeq \Z\x B\Ga^+_\infty$ which 
appeared in \cite{T1,T2}
relies on  a generalized group completion theorem. 
Tillmann constructs a homology fibration with fibre $\Z\x B\Ga_\infty$ 
and homotopy fibre of the homotopy type of $\Om B\ES$.  
The canonical map from the fibre to the canonical homotopy  
fibre is thus a homology equivalence.  
We use an explicit identification of the homotopy fibre with  
$\Om B\ES$ (before stabilization) 
to show that the  map given in Proposition \ref{equiS} 
induces the homology equivalence.

\vs

Consider the simplicial space 
with space of $n$-simplices  
$$(E_{\ES}\ES_1)_n=\coprod_{m_0,\dots,m_n\in Ob\ES} 
          \ES(m_0,m_1)\x\dots\x\ES(m_{n-1},m_n)\x\ES(m_n,1)$$ 
and boundary maps induced by composition in $\ES$. 
Consider also the telescope 
$$\ES_\infty(n)=Tel(\ES(n,1)\sta{T}{\rar}\ES(n,1)\sta{T}{\rar}\dots)$$
 where $\ES(n,1)\sta{T}{\rar}\ES(n,1)$ is induced by gluing the 
torus $T$. Note that 
$\ES_\infty(0)\simeq\Z\x B\Ga_\infty$. 
As composition induces maps 
$\ES(n,m)\x\ES_\infty(m)\to\ES_\infty(n)$, we can 
also  define a simplicial space  
$E_{\ES}\ES_\infty=Tel(E_\ES\ES_1\sta{T}{\rar}\dots)$. 
The map 
$$\pi:E_{\ES}\ES_\infty \rar B\ES,$$
induced by collapsing $\ES_\infty(n)$ to $\{n\}$, 
is a homology fibration.  
Let $hF_\infty:=PB\ES\x_{B\ES}E_{\ES}S_\infty$ denote the homotopy  
fibre. 
As $E_{\ES}\ES_\infty$ is contractible, $hF_\infty$  
is of the homotopy type of  
$\Om B\ES$. Hence, we have 
$$\Z\x B\Ga_\infty\simeq\ES_\infty(0)\sta{\simeq_{H_*}}{\rar} 
         hF_\infty\simeq\Om B\ES.$$

\begin{thm}\label{app}
The map $\psi:\ES(0,1)\to\Om B\ES$ defined by 
$\psi(\uF)=\uF\Big(^1_0\Big)D$ induces the homology equivalence 
$Z\x B\Ga_\infty\simeq_{H_*}\Om B\ES$, i.e. there  
is a commutative diagram
\[\xymatrix{\ES(0,1) \ar[rr]^\psi \ar[d] & & \Om B\ES\\
\ES_\infty(0) \ar[urr]_{\simeq_{H_*}} & &
}\]
\end{thm}

To prove the Theorem, we first study the non-stable case.  
Consider the map $\pi_{(1)}:E_{\ES}\ES_1\to B\ES$, and let  
$hF_1:=PB\ES\x_{B\ES}E_{\ES}\ES_1$ denote its homotopy fiber.  
 As $E_{\ES}\ES_1$ is also contractible, $hF_1$ is 
homotopy equivalent to $\Om B\ES$ (but not equivalent to 
the fiber in this case). 
 We want to construct an explicit  
homotopy equivalence  
$$\rho:hF_1\rar\Om_{01} B\ES$$
where $\Om_{01} B\ES$ is the space of loops in $B\ES$ starting at 
0 and ending at 1.
To a $q$-simplex  
$(\uF_1,\dots,\uF_q,\uG)\in   
   \ES(n_0,n_1)\x\dots\x\ES(n_{q-1},n_q)\x\ES(n_q,1)$
of $E_{\ES}\ES_1$ corresponds a $q+1$-simplex of $B\ES$
$$\s=n_0\sta{\uF_1}{\rar}n_1\rar\dots\sta{\uF_q}{\rar}n_q\sta{\uG}{\rar}1$$
having 1 as last vertex. The face opposite to 1 in $\s$ is 
$\pi_{(1)}\big(\uF_1,\dots,\uF_q,\uG\big)$. For $e\in E_{\ES}\ES_1$,  
define the path $\delta_e$ from $\pi(e)$ to 1,  
to be the straight line in $\s$ between $\pi(e)$ and 1 
(see figure). 
\begin{figure}[ht]
\begin{center}
\input line.pstex_t
\end{center}
\end{figure}
Now for $(p,e)\in hF_1$ ($e\in E_{\ES}\ES_1$ and $p$ is a path in  
$B\ES$ from 0 to $\pi(e)$),  
define  $\rho(p,e)$ to be the product of paths $p.\delta_e$.

\begin{lem}\label{hF1} 
The map $\rho:hF_1\rar\Om_{01} B\ES$ is a homotopy equivalence.
\end{lem}

\begin{proof}
Define the map  $\xi: \Om_{01}B\ES \rar hF_1$
by $\xi(\la)=(\la,Id_1)$, where $Id_1\in\ES(1,1)$ is the 
identity at 1. Then $\rho\circ\xi$
is the identity on $\Om_{01}B\ES$. On the other hand, $\xi\circ\rho$ is 
homotopic to the identity on $hF_1$. Indeed, for $e$ in the 
$q$-simplex $(\uF_1,\dots,\uF_q,\uG)$, we have 
$(\xi\circ\rho)(p,e)=(p.\de_e,Id_1)$. Consider the $q+1$-simplex 
of $E_{\ES}\ES_1$ defined by  $(\uF_1,\dots,\uF_q,\uG,Id_1)$. As in  
the case of $\de_e$, we can define a straight line $\ga_e$ 
in the $q+1$-simplex 
from $e$ to $Id_1$ ($e$ lies in the face opposite to $Id_1$). Now   
$\pi_{(1)}(\ga_e)=\de_e$. This induces the required
homotopy in $hF_1$ by truncating the path $p.\de_e$ at 
$\pi_{(1)}(\ga_e(t))$ (explicitly 
$H(t,p,e)=(p.\de|_{\pi_{(1)}(\ga_e(t))},\ga_e(t))$). 
\end{proof}

\noindent
{\em Proof of Theorem \ref{app}.} 
Consider the diagram 
\[\xymatrix{
\ES(0,1) \ar[r]^T \ar[d]_{\rho\circ j} 
              &\ES(0,1) \ar[r]^T \ar[d]_{\rho\circ j} &
               \ES(0,1) \ar[r] \ar[d]_{\rho\circ j} & \dots\\ 
\Om_{01}B\ES \ar[r]^{\la_T} & \Om_{01}B\ES \ar[r]^{\la_T} 
   &\Om_{01}B\ES \ar[r] & \dots
}\]
where $j:\ES(0,1)\to hF_1$ is the canonical map from the fiber to the 
homotopy fiber, $\rho$ and $T$ are defined above and 
$\la_T:\Om_{01}B\ES\to\Om_{01}B\ES$ is the multiplication with the 
loop from 1 to 1 defined by the torus $T$.  
As the squares commute only up to homotopy, we need to 
rectify this diagram to get a map of telescopes.

Let $\D$ be the discrete category with two copies of $\N$ as set of objects  
and morphisms as shown in the following diagram: 
\[\xymatrix{ A_0 \ar[r]\ar[d] & A_1 \ar[r]\ar[d] & A_2 \ar[r]\ar[d] 
                               & A_3 \ar[r]\ar[d]  & \dots\\
B_0 \ar[r] & B_1 \ar[r] & B_2 \ar[r]  & B_3 \ar[r]  & \dots
}\] 
where the squares commute strictly. 
Now consider the standard $k$-simplex 
$\Delta_k\subset\RR^{k+1}$ and consider the 
path space 
$$P_k=\{\delta:I\to\Delta_k|\delta(0)=(1,0,\dots,0)\ \  {\rm and} \ \ 
    \delta(1)=(0,\dots,0,1)\}.$$
So $P_k\simeq\Om\Delta_k$ is a contractible space.
Let $\tD$ be the category enriched over \Top \hspace{1mm}
having the same objects as 
$\D$ and morphism spaces defined as follows: for $n\in\N$ and $k\ge 0$,
$$\begin{tabular}{lcl}
 $\tD(A_n,A_{n+k})$ & = & $\D(A_n,A_{n+k})=\{*\}$\\
 $\tD(B_n,B_{n+k})$ & = & $\D(B_n,B_{n+k})=\{*\}$\\
 $\tD(A_n,B_{n+k})$ & = & $P_{k+1}$,
\end{tabular}$$
the other morphism spaces being empty. 
Labeling the $k+2$ vertices of $\Delta_{k+1}$ with $A_n,B_n,\dots,B_{n+k}$
induces a face inclusion $i:\Delta_{k+1}\hookrightarrow\Delta_{n+k+l+1}$.  
The composition of a morphism  
$f:A_n \to B_{n+k}$ with the unique morphism $g:B_{n+k}\to B_{n+k+l}$  
is defined to be the product of paths 
$(i\circ f).p$ in $\Delta_{n+k+l+1}$,   
where $p$ is the path from $B_{n+k}$ to 
$B_{n+k+l}$  following the edges 
$B_{n+k}\!-\!B_{n+k+1}\!-\dots-\!B_{n+k+l}$. 
To compose the  morphism $f:A_n \to A_{n+k}$ with a morphism 
$g:A_{n+k}\to B_{n+k+l}$, one uses the inclusion $j:\Delta_{l+1}\hookrightarrow 
\Delta_{k+l+1}$ sending $\Delta_{l+1}$ to the face  
having vertices labeled $A_n, B_{n+k},B_{n+k+1},\dots,B_{n+k+l}$. 
Define $g\circ f$ to be $j\circ g$.

The categories $\D$ and $\tD$ satisfy the hypothesis of Section \ref{diagsec}. 
We need to show that our data gives a functor $J:\tD\to\Top$. 
Define $J$ on object by 
$J(A_n):=\ES(0,1)$ and $J(B_n):=\Om_{01}B\ES$.
The diagram given at the beginning of the proof defines $J$ on the 
morphisms of the type $A_n\to A_{n+k}$, $B_n\to B_{n+k}$ and $A_n\to B_n$  
as well as their composition. Such  
compositions between $A_n$ and $B_{n+k}$ are  paths following the 
edges in the relevant simplex. The figure shows in the case  
of $A_n\to B_{n+2}$ that the images are actually paths following the  
edges in a simplex of $B\ES$. 
\begin{figure}[ht]
\begin{center}
\input tele.pstex_t
\end{center}
\end{figure}
For any morphism $\de:A_n \rar B_{n+k}$, define
$J(\de): \ES(0,1)\rar \Om_{01}B\ES$
by setting $J(\de)(\uF)$ to be the corresponding path in the 
$k+1$-simplex of $B\ES$
$$0\sta{\uF}{\rar}1\sta{T}{\rar}1\sta{T}{\rar}\dots\sta{T}{\rar}1.$$
$J$ is functorial essentially because we defined composition in 
$\tD$ precisely to make it functorial.

Let $J'=p_*J:\D\rar \Top\ $ 
be the rectification of $J$.  
The rectification produces two telescopes $J'(A_*)$ and $J'(B_*)$ 
which are equivalent  
to the one we started with by naturality of the equivalence  
$p^*p_*J\simeq J$. Moreover, we   
now have a map of telescopes $f:Tel(J'(A_*))\to Tel(J'(B_*))$. 
From \cite{T1}, we know that the map 
$ \ES_\infty(0)\to hF_\infty=PB\ES\x_{B\ES}E_{ \ES} \ES_\infty$
is a homology equivalence. Putting all this information together, we have
\[\xymatrix{Tel(J'(A_*)) \ar[d] \ar@{<->}[r]^\simeq 
           & \ES_\infty(0) \ar[r]^{\simeq_{H_*}} 
           & hF_\infty\\
Tel(J'(B_*)) \ar@{<->}[r]^\simeq 
           & Tel(\Om_{01}B\ES) &   
}\]
Now on each ``level'' of the telescope, by naturality of the equivalence and 
by Lemma \ref{hF1}, we have a commutative diagram
\[\xymatrix{J'(A_n) \ar@{<->}[r] \ar[d] & \ES(0,1)\ar[d]_{\rho\circ j}\ar[r]^j 
          & hF_1 \ar[dl]^\rho_{\simeq}\\
J'(B_n) \ar@{<->}[r] & \Om_{01}B\ES &
}\]
It follows  that the map 
$f:Tel(J'(A_*))\to Tel(J'(B_*))$ is a homology 
equivalence.  
Hence we have a commutative diagram
\[\xymatrix{
& \ES(0,1) \ar[r] \ar[d]^{\rho\circ j}\ar[dl]_\psi  
                        & \ES_\infty(0) \ar[d]^{\simeq_{H_*}}\\  
\Om B\ES & \ar[l]_\simeq \Om_{01}B\ES \ar[r]^{\simeq\ \ \ }  
      & Tel(\Om_{01}B\ES) }\]
where the map
$\Om_{01}B\ES\to\Om B\ES$ is the multiplication with the path  
from 0 to 1 along the disc (taken backwards), and 
$Tel(\Om_{01}B\ES)\simeq\Om_{01}B\ES$ as the telescope structure 
map has a homotopy inverse. 
\cqfd

}

\end{document}